\documentclass{amsart}
\usepackage[all]{xy}
\numberwithin{equation}{section}
\usepackage{amssymb}
\let\cal\mathcal

\def\Nscr{{\cal N}}
\def\Oscr{{\cal O}}
\def\Pscr{{\cal P}}

\let\blb\mathbb

\def\EE{{\blb E}}

\def \ZZ{{\blb Z}}

\def \NN{{\blb N}}

\def \HH{{\blb H}}

\def\Lotimes{\overset{L}{\otimes}}

\def\gr{\operatorname{gr}}

\def\coh{\mathop{\text{\upshape{coh}}}}

\def\gr{\operatorname {gr}}
\def\Spec{\operatorname {Spec}}

\def\Ext{\operatorname {Ext}}
\def\Hom{\operatorname {Hom}}

\def\End{\operatorname {End}}
\def\RHom{\operatorname {RHom}}
\def\uRHom{\operatorname {R\mathcal{H}\mathit{om}}}
\def\Sl{\operatorname {Sl}}

\def\End{\operatorname {End}}

\def\rk{\operatorname {rk}}

\def\gldim{\operatorname {gl\,dim}}

\def\r{\rightarrow}

\newtheorem{lemma}{Lemma}[section]
\newtheorem{proposition}[lemma]{Proposition}
\newtheorem{theorem}[lemma]{Theorem}
\newtheorem{corollary}[lemma]{Corollary}

\newtheorem{lemmas}{Lemma}[subsection]
\newtheorem{propositions}[lemmas]{Proposition}
\newtheorem{theorems}[lemmas]{Theorem}
\newtheorem{corollarys}[lemmas]{Corollary}

\theoremstyle{definition}

\newtheorem{examples}[lemmas]{Example}

{

}

\theoremstyle{remark}

\newdimen\uboxsep \uboxsep=1ex
\def\uboxn#1{\vtop to 0pt{\hrule height 0pt depth 0pt\vskip\uboxsep
\hbox to 0pt{\hss #1\hss}\vss}}

\def\uboxs#1{\vbox to 0pt{\vss\hbox to 0pt{\hss #1\hss}
\vskip\uboxsep\hrule height 0pt depth 0pt}}

\def\SL{\operatorname{SL}}
\def\red{\operatorname{red}}
\title{Some new examples of non-degenerate quiver potentials}
\author{L. de Thanhoffer de V\"olcsey}
\author{M. Van den Bergh}
\email{louicious@mac.com, michel.vandenbergh@uhasselt.be}
\address{Departement WNI,  Universiteit Hasselt,
 3590 Diepenbeek, Belgium.}
\thanks{The second author is a senior researcher at the FWO}
\begin{document}
\begin{abstract}
  We prove a technical result which allows us to establish the
  non-degeneracy of potentials on quivers in some previously unknown
  or non-obvious cases.

Our result applies to certain McKay quivers and also to potentials
derived from geometric helices on Del Pezzo surfaces. 
\end{abstract}
\maketitle
\section{Introduction}
In \cite{DWZ} Derksen, Weyman and Zelevinsky describe how to
``mutate'' at a vertex a pair $(Q,W)$ consisting of a quiver $Q$
and a potential $W\in (kQ)\,\hat{}/[(kQ)\,\hat{},(kQ)\,\hat{}]$. This construction
produces a new such pair $(Q',W')$. The Jacobian (or Ginzburg)
algebras of the pairs $(Q,W)$, $(Q',W')$ share many homological
properties \cite{BIRS,IRe,KY}.

One peculiarity of the mutation process is that
it is only defined if the vertex is not incident to a loop or two-cycle.
Even if  all vertices in $Q$ have this property then this  is not
necessarily the case for $Q'$. 
If the property of having no loops or two-cycles persists under iterated mutations
then we say that $(Q,W)$ (or $W$) is non-degenerate.

In this paper we give a technical result (Theorem \ref{ref-3.1-3})
which allows us to establish the non-degeneracy of potentials in some
previously unknown or non-obvious cases. 
Our theorem implies for example that potentials derived from geometric helices
on Del Pezzo surfaces (see below) are non-degenerate. 
In this way we recover part of the
main result of \cite{BS1}.

\medskip

Our
theorem also applies  to the McKay quiver associated to a group
$G\subset\SL(V)$ where $\dim V=3$ and $G$ acts freely on
$V^\ast-\{0\}$ (this implies that $G$ is cyclic).  See Corollary
\ref{ref-4.1.2-8} below.

In fact the original motivation of our work was an  explicit example considered by Iyama and Reiten in \cite{IRe}
where they mention that they do not know if the 
mutations of the associated McKay quiver do not create loops or two-cycles.

Assume that $Q$ is
the McKay quiver of the pair $(\ZZ/5\ZZ,kx_1+kx_2+kx_3)$ where the
generator of $\ZZ/5\ZZ$ acts diagonally by $(\xi,\xi^2,\xi^2)$. The
quiver $Q$ looks as follows
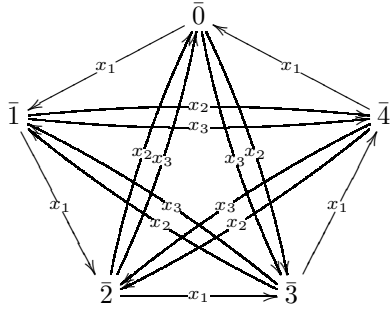
\[
\xymatrix{
&&\bar{0}\ar[lld]|{x_1}\ar@/^0.4em/[dddr]|{x_2}\ar@/_0.4em/[dddr]|{x_3}&&\\
\bar{1}\ar@/^0.4em/[rrrr]|{x_2}\ar@/_0.4em/[rrrr]|{x_3}\ar[rdd]|{x_1}&&&&
\bar{4}\ar@/^0.4em/[ddlll]|(0.55){x_2}\ar@/_0.4em/[ddlll]|(0.55){x_3}
\ar[llu]|{x_1}\\
&&&&\\
&\bar{2}\ar@/^0.4em/[ruuu]|{x_2}\ar@/_0.4em/[ruuu]|{x_3}\ar[rr]|{x_1}&&
\bar{3}\ar[ruu]|{x_1}\ar@/^0.4em/[llluu]|(0.45){x_2}\ar@/_0.4em/[llluu]|(0.45){x_3}&
}
\]
The potential $W$ is the signed sum of all three-cycles containing
$x_1,x_2$ and $x_3$.  The Jacobian algebra of $(Q,W)$ is the skew
group ring $k[[x_1,x_2,x_3]]\#\ZZ/5\ZZ$.  Our main theorem implies that
$W$ is non-degenerate.

\medskip

On the other hand if $G$ does not act freely on $V^\ast-\{0\}$ then
the associated potential will generally be degenerate. See Example \ref{ref-4.1.3-9}
below.

\medskip

The main technical tool in this paper is the lifting of mutations
to the graded setting (see \cite{AO}).

\medskip

Throughout all quivers are finite. $k$ is a base field.
\section{Graded and ungraded mutations}
Let $Q$ be a quiver. A potential on $Q$ is an element $W\in
(kQ)\,\hat{}/[(kQ)\,\hat{},(kQ)\,\hat{}]$ containing no paths of length one
or zero. We write $W$ as a sum of oriented cycles in $Q$. The
corresponding Jacobian algebra is defined in the usual way
\[
\Pscr(Q,W)=(kQ)\,\hat{}/(\partial_a W)_a
\]
The space of paths
of length zero in $kQ$ will be denoted by $R$. This is a semisimple $k$-algebra.

A potential is said to be reduced if it contains no two-cycles. Such
two-cycles lead to relations which express some arrows in terms of
others and hence such arrows can be eliminated in the Jacobian algebra. This
observation is refined in the splitting theorem \cite[Thm
4.6]{DWZ} which asserts that $(Q,W)$ is ``right equivalent'' (see \cite[Def.\
  4.2]{DWZ}) to a direct sum decomposition
\begin{equation}
\label{ref-2.1-0}
(Q,W)\sim (Q^{\text{red}},W^{\text{red}})\oplus (Q^{\text{triv}},W^{\text{triv}})
\end{equation}
where $W^{\text{red}}$ is reduced and $W^{\text{triv}}$ contains only
two-cycles and its associated Jacobian algebra is equal to $R$.

The decomposition theorem implies in particular 
\begin{equation}
\label{ref-2.2-1}
\Pscr(Q,W)\cong\Pscr(Q^{\text{red}},W^{\text{red}})
\end{equation}
If $A$ is an abelian group then we say that $Q$ is $A$-graded if we
have assigned an $A$-degree $|a|$ to all arrows $a$ in $Q$. If $W$ is
homogeneous of degree $r\in A$ the partially completed Jacobian algebra
\[
\Pscr^{\gr}(Q,W)=(kQ)^{\gr}\,\hat{}/(\partial_a W)_a
\]
(where we complete only with respect to sequences of paths having ascending length
but constant degree),
is naturally an $A$-graded
algebra. It is observed in \cite{AO} that the graded analogue of the splitting
theorem holds and hence the decomposition \eqref{ref-2.1-0} can be 
performed on the graded level.

The following lemma is very useful.
\begin{lemma} 
\label{ref-2.1-2}  Let $Q'$ be obtained from $Q$ by repetively deleting pairs
of arrows $i\xrightarrow{a} j\xrightarrow{b} i$ with $i\neq j$ and $|a|+|b|=r$.
If $Q'$ contains a two-cycle then so does $Q^{\red}$. 
\end{lemma}
\begin{proof}
By the degree constraint on $W$ we can only eliminate two-cycles of the form $ab$ with
$|a|+|b|=r$. This implies the result. 
\end{proof}

\medskip

Assume that $Q$ does not have loops or two-cycles and let $i$ be a
vertex of $Q$. The mutation $\mu_i(Q,W)=(Q',W')$ at $i$ is
defined in \cite{DWZ} as follows.
\begin{itemize}
\item For any sequence of arrows $u\xrightarrow{a}i\xrightarrow{b}v$ we add an
arrow $u\xrightarrow{[ba]}v$.
\item All arrows $a$ starting or ending in $i$ are replaced with
  opposite arrows $a^\ast$.
\end{itemize}
We have
\[
W'=[W]+\sum b^\ast a^\ast[ab]
\]
where $[W]$ is obtained from $W$ by replacing all compositions $ba$
through the vertex $i$ by the new arrows $[ba]$.
Following \cite{DWZ} we put $\tilde{\mu}_i(Q,W)=(Q^{\prime \red},W^{\prime\red})$. 
The operations $\tilde{\mu}_i$ is an operation on right equivalence
classes of quivers with potential (see \cite[Thm 5.2]{DWZ}). 

Now assume that $Q$ is in addition $A$-graded and that $W$ is homogeneous of degree
$r$. Then we make $Q'$ into a $A$-graded quiver by fixing the degrees of the arrows
as follows \cite{AO}
\begin{itemize}
\item Arrows which are both in $Q$ and $Q'$ do not change degrees.
\item $|[ab]|=|a|+|b|$.
\item $|a^\ast|=-|a|+r$ if $a$ ends in $i$.
\item $|a^\ast|=-|a|$ if $a$ starts in $i$. 
\end{itemize}
With this grading $W'$ is homogeneous of degree $r$.  If we perform
mutation on the graded level we write $\mu_i^{\gr}(Q,W)=(Q',W')$
and $\tilde{\mu}_i^{\gr}(Q,W)=(Q^{\prime\red},W^{\prime\red})$. Such graded
mutations are compatible with forgetting the grading. 
\section{Main theorem}
Let $S=\bigoplus_{n\ge 0} S_n$ be a commutative noetherian $\ZZ$-graded ring and let $\Lambda$
be a graded $S$-order (i.e.\ $\Lambda$ is a graded $S$-algebra, finite as a module over
$S$). We will say that $\Lambda/S$ is almost Azumaya
if for all $P\in \Spec S$ such that $P\not\supset S_{\ge 1}$ the algebra
$\Lambda_P$ is Azumaya over $S_P$. This implies
in particular that for any non-zero idempotent $e\in \Lambda_0$ the quotient $\Lambda/\Lambda e\Lambda$
is a finite $S_0$-module. The latter property is in fact the only one we will use below. 
\begin{theorem}
\label{ref-3.1-3}
Assume the following assumptions hold
\begin{enumerate}
\item
$(Q,W)$ is a $\ZZ$-graded connected quiver with a reduced homogeneous
potential~$W$ of degree $r$.
\item $Q$ has at least three vertices.
\item Put $\Lambda=\Pscr^{\gr}(Q,W)$. Then $\dim \Lambda_i<\infty$ for all $i$ and
$\dim \Lambda_i=0$ for $i\ll 0$. 
\item $\Lambda$ is  a 3-dimensional CY-order over a noetherian center which is almost
  Azumaya.
\item The zeroeth Hochschild homology $\operatorname{HH}_0(\Lambda)\overset{\text{def}}{=}\Lambda/[\Lambda,\Lambda]$ of $\Lambda$
  contains no homogeneous elements with degree lying in the interval
  $[1,r/2]$.
\end{enumerate}
Then $Q$ has no loops or two-cycles.
\end{theorem}
\begin{proof}
  Assume first that $Q$ contains a loop $i\xrightarrow{a} i$ and let
  $e$ be the sum of the idempotents corresponding to the vertices
  different from $i$. Let $Q_0$ be the subquiver of $Q$
  consisting of all loops at $i$. Put $\bar{\Gamma}=\Gamma/\Gamma
  e\Gamma$. Then $\bar{\Gamma}$ is the graded Jacobian algebra
  $\Pscr^{\text{gr}}(Q_0,W_0)$ where $W_0$ is obtained from $\bar{W}$
  by dropping all cycles passing through vertices other than $i$.

If $W_0=0$ then $\bar{\Gamma}$ is infinite dimensional which is impossible
by the discussion preceding the theorem. Hence $W_0$ is a sum of terms of
degree $r$ which
are products of loops
\begin{equation}
\label{ref-3.1-4}
a_1\cdots a_n
\end{equation}
with $n\ge 3$ (since $W_0$
is reduced). Assume that \eqref{ref-3.1-4} is the shortest such term. Then there must
be some $a_i$ which has degree $\le r/3\le r/2$. As the paths occurring in
the relations of $\Gamma/\Gamma e \Gamma$ are products of at least two arrows we
have $a_i\neq 0$ in $\Gamma/\Gamma e \Gamma$.  Furthermore it is also clear that $a_i$ is
not in the image of $[kQ_0,kQ_0]$. Thus $a_i$ represents a non-zero
element of $\operatorname{HH}_0(\bar{\Gamma})$. Since there is a surjection
$\operatorname{HH}_0(\Gamma)\r \operatorname{HH}_0(\bar{\Gamma})$, we
obtain a non-trivial element of $\operatorname{HH}_0(\Gamma)$ as well. 
Since $|a_i|\le r/2$ we obtain from the hypotheses that $|a_i|\le 0$. But
then $a_i$ must be nilpotent in $\Gamma$ and hence by \cite{LH} $a_i$
is a sum of commutators, which is a contradiction.

We now assume that $Q$ has no loops. The proof in the case of a
two-cycle $i\rightarrow j\rightarrow i$ is similar and we only sketch it.  Let $Q_0$ be the
subquiver of $Q$ consisting of all arrows between $i$ and $j$
and vice versa and let $e$ be the sum of the idempotents corresponding
to the vertices different from $i$ or $j$.  Put
$\bar{\Gamma}=\Gamma/\Gamma e\Gamma$.  Then
$\bar{\Gamma}=\Pscr^{\text{gr}}(Q_0,W_0)$ where $W_0$ is obtained from
$W$ by dropping all cycles passing through vertices other than
$i$ or $j$.

The case $W_0=0$ is again impossible. Hence $W_0$ is a sum of terms of
degree $r$ which
are products of two-cycles
\[
a_1b_1a_2b_2\cdots a_nb_n
\]
with $i\xrightarrow{a_l} j$ and $j\xrightarrow{b_m} i$ with $n\ge 2$. Hence
one of the two-cycles $a_ib_i$ must have degree zero. Since this time
the relations have degree $\ge 3$ we have $a_ib_i\neq 0$. This gives
a non-trivial element in $\operatorname{HH}_0(\Gamma)$ which must have degree $\le 0$ by
the hypotheses. Thus
$a_ib_i$ is nilpotent and hence is a sum of commutators which is a contradiction.

\end{proof}
\begin{proposition} \label{ref-3.2-5} Let $(Q,W)$ be as in Theorem \ref{ref-3.1-3} and let 
$i$ be a vertex of $Q$. Then up to right equivalence we may assume
that $\tilde{\mu}_i(Q,W)=(Q',W')$ also satisfies the conditions of
Theorem \ref{ref-3.1-3}. In particular $Q'$ has no loops or two-cycles. 
\end{proposition}
\begin{proof}
Since ungraded and graded mutations are
  compatible and since the decomposition \eqref{ref-2.1-0} can be
  performed at the graded level, we can up to right equivalence,
  assume that $Q'$ is graded and  $W'$ is homogeneous.

  It has been proved in \cite[Thm 5.2]{BIRS} and in
  \cite[\S6]{KY} that for Jacobian algebras which are
  Calabi-Yau of dimension three, Derksen-Weyman-Zelevinsky mutations
  correspond to Iyama-Reiten mutations. With some work
  these proofs can be carried over to the graded context.  Hence in
  particular we obtain that
  $\Gamma=\Pscr^{\text{gr}}(Q',W')$ is 3-CY. Furthermore
  $Z(\Gamma)=Z(\Lambda)$, $\Gamma$ is almost Azumaya over its center and $\Gamma$ is derived equivalent to
  $\Lambda$ \cite{IRe}.  Finally by the derived invariance of
  Hochschild homology we have an isomorphism of graded vector spaces
\[
\Lambda/[\Lambda,\Lambda]=\operatorname{HH}_0(\Lambda)=\operatorname{HH}_0(\Gamma)=
\Gamma/[\Gamma,\Gamma]\qed
\]
\def\qed{}\end{proof}
\begin{corollary}
\label{ref-3.3-6}
If $(Q,W)$ is as in Theorem \ref{ref-3.1-3} then $(Q,W)$ is non-degenerate.
\end{corollary}
\begin{proof} This is clear from Proposition \ref{ref-3.2-5}.
\end{proof}
\section{Examples and counterexamples}
\subsection{Skew group rings}
In this section we assume that $k$ is algebraically closed of
characteristic zero. Let $V$ be a three dimensional vector space and
let $G$ be a finite subgroup of $\SL(V)$ of order $n$. Then it is
well-known that the skew group algebra $\Lambda=SV\# G$ is
$3$-Calabi-Yau. Furthermore $\Lambda$ is graded Morita equivalent to a
graded Jacobian algebra $\Pscr^{\text{gr}}(Q,W)$ where $Q$ is the
McKay quiver of $(G,V)$ and $W$ has degree three. The McKay quiver is
the quiver with vertices the irreducible representations $(V_i)_i$ of
$G$ and arrows $i\r j$ with multiplicity $p$ if $V_j$ occurs
$p$ times in $V\otimes V_i$. 
We refer to \cite[Thm 3.2]{BSW} for a detailed discussion
on how to construct the potential.


In the case of a cyclic group $G\cong \ZZ/n\ZZ$, the situation is
particularly simple.  The irreducible representations of $G$ are given
by characters (``weights''). After the choice of a primitive $n$-th root of
unity $\xi$ a character $\chi$ may be identified with an element $\bar{a}$ of
$\ZZ/n\ZZ$ through the rule $\chi_{\bar{a}}(\bar{m})=\xi^{am}$.

 Let $\bar{a}_1$, $\bar{a}_2$, $\bar{a}_3$ be the weights of $V$.  The
 fact that $G\subset \SL(V)$ is equivalent to $\sum_i
 \bar{a}_i=\bar{0}$.  The vertices of the McKay quiver of $(G,V)$
 are indexed by the elements of $\ZZ/n\ZZ$. If $\bar{l}\in \ZZ/n\ZZ$
 then there are three arrows $x_1,x_2,x_3$ leaving $\bar l$ and ending
 respectively in $\bar{l}+\bar{a}_1$, $\bar{l}+\bar{a}_2$,
 $\bar{l}+\bar{a}_3$. The potential is given by
\begin{equation}
\label{ref-4.1-7}
W=\sum \pm x_{i_1}x_{i_2}x_{i_3}
\end{equation}
where the sum runs over all three-cycles such that
$\{{i_1},{i_2},{i_3}\}=\{1,2,3\}$ and the sign is positive iff
$(i_1,i_2,i_3)$ is equal to $(1,2,3)$, up to cyclic permutation.

\begin{lemmas}
The hypotheses of Theorem \ref{ref-3.1-3} hold for $(Q,W)$ if and only if $G\cong \ZZ/n\ZZ$
is cyclic and acts with weights $(\bar{a}_i)_i$ on $V$ such that $\gcd(n,a_i)=1$.
\end{lemmas}
\begin{proof}
  Assume that the hypotheses for Theorem \ref{ref-3.1-3} hold.  It is
  well-known that $\Lambda$ is almost Azumaya over its center if and
  only if $SV^G$ has an isolated singularity which is equivalent to $G$
  acting freely on $V^\ast-\{0\}$ and
hence that $G$ is
  cyclic (see \cite{KN}). The fact that $G$ acts freely is equivalent to the stated condition
on weights.

Conversely assume that $G$ is cyclic and acts with weights relatively prime to~$n$. We
verify the hypotheses of Theorem \ref{ref-3.1-3}: we have already mentioned (1).
(2) follows from the fact that $\ZZ/2\ZZ$ cannot act freely
 on $V-\{0\}$ (one weight must be zero). (3) is
clear. (4) has been mentioned previously.
It remains to check (5). We must show $\operatorname{HH}_0(\Lambda)_1=0$.
If $\operatorname{HH}_0(\Lambda)_1\neq 0$ then $Q$ must have a loop. This implies that one
of the weights is zero and thus not relatively prime to $n$. Contradiction!
\end{proof}
We obtain:
\begin{corollarys} \label{ref-4.1.2-8}
Assume that $Q$ is the McKay quiver of a cyclic group $\ZZ/n\ZZ$
acting with weights $\bar{a}_1,\bar{a}_2,\bar{a}_3$ such that $\gcd(a_i,n)=1$ and
$\sum_i\bar{a}_i=0$. Let $W$ be the potential \eqref{ref-4.1-7}. Then
$W$ is non-degenerate.
\end{corollarys}
%
On the other hand the hypotheses for Corollary \ref{ref-4.1.2-8} can not be weakened
as shown by the following example:
\begin{examples} \label{ref-4.1.3-9} Let $G=\ZZ/6\ZZ$ and assume that $G$ acts with
  weights $(\bar{2},\bar{5},\bar{5})$. Since $2{\mid} 6$ the
  hypotheses for Corollary \ref{ref-4.1.2-8} do not hold.

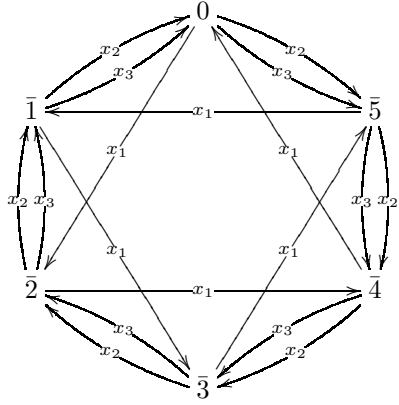
\[
\xymatrix{
&&\bar{0}\ar[dddll]|{x_1}\ar@/^0.5em/[drr]|{x_2}\ar@/_0.5em/[drr]|{x_3}\\
\bar{1}\ar[rrddd]|{x_1}\ar@/^0.5em/[urr]|{x_2}\ar@/_0.5em/[urr]|{x_3} &&&&
\bar{5}\ar[llll]|{x_1}\ar@/^0.5em/[dd]|{x_2}\ar@/_0.5em/[dd]|{x_3}\\
&&&&&\\
\bar{2}\ar[rrrr]|{x_1}\ar@/^0.5em/[uu]|{x_2}\ar@/_0.5em/[uu]|{x_3}&&&&
\bar{4}\ar[lluuu]|{x_1}\ar@/^0.5em/[dll]|{x_2}\ar@/_0.5em/[dll]|{x_3}\\
&&
\bar{3}\ar[rruuu]|{x_1}\ar@/^0.5em/[ull]|{x_2}\ar@/_0.5em/[ull]|{x_3}&&
}
\]
We note that the McKay quiver of a cyclic groups is naturally
$\ZZ^3$-graded. It is convenient to use the grading by the monomials
in $x_1,x_2,x_3$ such that $|x_i|=x_i$. Then $W$ is homogeneous of
degree $x_1x_2x_3$. With these conventions we see that performing a mutation
at the vertex $\bar{0}$ results in a two-cycle $[x_1x_1]x_1$ between vertices $\bar{2}$ and
$\bar{4}$. This two-cycle has degree $x_1^3$ and cannot be eliminated
for degree reasons (see Lemma \ref{ref-2.1-2}). 
\end{examples}
\subsection{Del Pezzo surfaces}
\label{ref-4.2-10}
Let $Y$ be a Del-Pezzo surface and let $(E_i)_{i=1,\ldots,n}$ be a
full exceptional collection on~$Y$. Put $\EE=\bigoplus_{i=1}^n E_i$
and define $A(\EE)=\End(\EE)$. Assume that
\[
\HH(E)=(E_i)_{i\in\ZZ}=(\ldots,\omega_Y\otimes E_n,E_1,\ldots,E_n,\omega^{-1}_Y\otimes, E_1,\ldots)
\]
is a geometric helix (see \cite{Bondal,BS1}). I.e.\ every slice
$(E_{i+1},\ldots,E_{i+n})$ is an exceptional collection and
furthermore
\[
\forall i<j, \forall k>0:\Ext^k(E_i,E_j)=0
\]
The rolled up helix algebra is  the $\ZZ$-graded ring
\[
B(\HH)=\bigoplus_{k\in \ZZ} \Hom(\EE,\omega_Y^{-k}\otimes\EE)
\]
with obvious multiplication. 
\begin{theorems} \label{ref-4.2.1-11} $B(\HH)$ is an $\NN$-graded Jacobian algebra derived
  from a graded super potential of degree one which is Calabi-Yau of
  dimension three (as ungraded algebra) and which is almost Azumaya over its center.
\end{theorems}
This result is folklore among physicists. It follows from 
the work by Ed Segal in \cite{segal} although it is not really stated
explicitly there. The fact that $B(\HH)$ is Calabi-Yau of dimension dimension three is
\cite[Thm 3.6]{BS1}. 
 The fact that $B(\HH)$ is Jacobian follows from \cite{VdBsuper} and probably also
from the proof of \cite[Thm 3.1]{Bocklandt} although the formulation
this theorem is for quivers whose arrows have degree one. 

We give for the convenience of the reader an
independent proof of this theorem in Appendix~\ref{ref-A-15} based on some results from
\cite{Keller11}.

\begin{theorems} \label{ref-4.2.2-12} Let $(Q,W)$ be such that $W$ is reduced
  and $\Pscr^{\gr} (Q,W)=B(\HH)$.  Then $(Q,W)$ satisfies the
  hypotheses of Theorem \ref{ref-3.1-3}.
\end{theorems}
In particular it follows from Corollary \ref{ref-3.3-6} 
 that $(Q,W)$ is non-degenerate. This is a consequence of \cite[Thm 1.7]{BS1}. In loc.\ cit.
the authors deduce this fact from their result that a mutation of $B(\HH)$ is always
of the form $B(\HH')$ for another geometric helix $\HH'$.  
\begin{proof}[Proof of Theorem \ref{ref-4.2.2-12}] (1) is true with $r=1$. (2) is true since $\rk K_0(\coh(Y))\ge 3$.
(3) is obvious. (4) follows from Theorem \ref{ref-4.2.1-11}. (5) is vacuous
since the potential has degree one. 
\end{proof}
\subsection{Deformed preprojective algebras}
In this section we given an easy example which shows that hypothesis
(5) of Theorem \ref{ref-3.1-3} cannot be omitted. For simplicity we
assume that $k$ is algebraically closed of characteristic zero.

Let $Q=(Q_0,Q_1)$ ($Q_0=$vertices, $Q_1=$edges) be an extended Dynkin
quiver and let $\tilde{Q}$ be its associated double quiver. That is,
for every arrow $a\in Q_1$ we adjoin an arrow $a^\ast$ with the opposite
orientation. By virtue of the construction every arrow in $\tilde{Q}$ is
part of a two-cycle.

For $(\lambda_i)_i\in
k^{Q_0}$ we define a corresponding central extension of the
preprojective algebra as follows
\[
\Pi^{\lambda}=k[t]\tilde{Q}\left/\biggl(\sum_{a\in Q_1}[a,a^\ast]-\sum_{i\in Q_0} \lambda_i t e_i\biggr)\right.
\]
where as usual $e_i$ is the idempotent corresponding to $i\in
Q_0$. This algebra (or at least its principal representation space) 
is introduced in \cite{CKV}. It has a natural grading
with the arrows in $\tilde{Q}$ having degree one and $\deg t=2$.  Let $\bar{Q}$ be
obtained from $\tilde{Q}$ by adjoining an additional loop at every
vertex. Then $\Pi^{\lambda}$ is the graded Jacobian algebra of
$(\bar{Q},\bar{W})$ with
\[
\bar{W}=\left(\sum_{i\in Q_0}t_i\right)\cdot\left(\sum_{a\in Q_1}[a,a^\ast]\right)-\sum_{i\in Q_0} \frac{1}{2}\lambda_i t_i^2 
\] 
where $t_i=te_i$ is the added loop at $i$. This induces the potential 
\cite[eq.\ (3.3)]{CKV} on representation spaces. 

If all $\lambda_i$ are non-zero then the loops can be eliminated and $\Pi^\lambda$
becomes the Jacobian algebra of $(\tilde{Q},W)$ with 
\[
W=\sum_{a,b\in Q_1,i\in Q_0} \frac{1}{\lambda_i} e_i[a,a^\ast][b,b^\ast]
\]
Let $0$ be the extending vertex of $Q$ (i.e.\ $Q-\{0\}$ is a Dynkin diagram).
Following \cite{CH} we  say that $\alpha\in \NN Q_0$ is a Dynkin root
if $\alpha$ is a root such that $\alpha_0=0$.

The following proposition collects some basic facts about $\Pi^\lambda$.
\begin{propositions} 
\label{ref-4.3.1-13} Let $\delta\in kQ_0$ be the
  indecomposable imaginary root and assume that $\delta\cdot
  \lambda\overset{\text{def}}{=}\sum_{i\in Q_0}\delta_i\lambda_i=0$.
  Assume in addition that for all Dynkin roots $\alpha$ one has
  $\alpha\cdot\lambda\neq 0$ (this is implies $\lambda_i\neq 0$ for all $i$). Then $\Pi^\lambda$ is a 3-CY order 
  which is almost Azumaya over a noetherian center.
\end{propositions}
\begin{proof} Most of this is stated in \cite{QB} and proved
  using results from \cite{CH}. 
Let $G\subset \Sl_2(k)$ be the Kleinian group whose McKay quiver is $\tilde{Q}$.
The vertices of $\tilde{Q}$ correspond to the irreducible
representations of $G$ and in particular the extending vertex corresponds
to the trivial representation.

For $i\in Q_0$ let $f_i$ be the central idempotent in $kG$
corresponding to the irreducible $G$-representation indexed by $i$. 
Then
$\Pi^\lambda$ is Morita equivalent to another ring $A^{\lambda}$ (a
symplectic reflection algebra)
\[
k[t]\langle x,y\rangle\#G/\biggl(xy-yx-t\sum_{i\in Q} \frac{\lambda_i}{\delta_i} f_i\biggr)
\]
and under this equivalence the augmentation idempotent $e\in kG$ corresponds to $e_0$
(see \cite[Thm 3.4]{CH}). The condition $\delta\cdot\lambda=0$ translates
into the fact that $\sum_{i\in Q} \frac{\lambda_i}{\delta_i} f_i$ has
trace zero for the regular representation. In other words its image
under the projection on the identity component $kG\r k$ is zero.

Based on the fact that $\gr A^{\lambda}=k[t][x,y]\#G$ it is deduced in
\cite[\S2.2]{QB} that $A^{\lambda}$ is a $k[t]$-flat prime noetherian
maximal order which is Auslander-Gorenstein and Cohen-Macaulay of GK
 and global dimension three. In a similar way $C^{\lambda}$ is
a noetherian integral domain of GK dimension 3 which is a flat
(graded) deformation of the hypersurface ring $k[x,y]^G$.  In
particular $C^{\lambda}$ is also Gorenstein. Furthermore one has
\begin{equation}
\label{ref-4.2-14}
A^{\lambda}\cong \End_{C^{\lambda}}(A^{\lambda}e)
\end{equation}
In \cite[Cor.\ 3.5]{QB} it is shown that $C^{\lambda}$ is
commutative and that $C^\lambda$ is in fact isomorphic to the center of
$A^{\lambda}$. Finally it then follows from \cite[Prop.\ 3.2]{QB}
that $A^{\lambda}$ is a non-commutative crepant resolution of its center \cite{VdB32}. 
This implies that $A^{\lambda}$ is 3-CY.  All these results carry over to
$\Pi^\lambda$ by Morita equivalence.

In particular \eqref{ref-4.2-14} becomes 
\[
\Pi^{\lambda}\cong \End_{e_0\Pi^\lambda e_0}(\Pi^{\lambda}e_0)
\]
with $e_0\Pi^\lambda e_0\cong C^{\lambda}$.  To prove that $\Pi^\lambda$
is almost Azumaya it suffices to prove that
$S=\Pi^{\lambda}/\Pi^{\lambda}e_0\Pi^{\lambda}$ is finite dimensional.
According to \cite[Cor.\ 9.6]{CH} $S/S(t-1)=0$. Since $S$ is graded
this is implies that $S_t=0$ (this is a statement purely about
graded $k[t]$-modules). Since $S$ is finitely generated as $\Pi^\lambda$-module it
follows that $t^n S=0$ for some $n$. On the other hand we also have that
$S/tS$ is finite dimensional since for example $\Pi^\lambda/t\Pi^\lambda$ is
the preprojective algebra of $Q$ which is Morita equivalent to $k[x,y]\#G$
and the latter is almost Azumaya (see \cite{ReVdB}). It now follows
that $S$ is finite dimensional.
\end{proof}
So if $Q$ has three or more vertices (i.e.\ $Q\neq A_1,A_2$) and $\lambda$ is as in
Proposition \ref{ref-4.3.1-13} then $\Pi^\lambda$ satisfies all conditions of
Theorem \ref{ref-3.1-3} except (5). Indeed the degree of the potential
is four and $\Pi^\lambda/[\Pi^\lambda,\Pi^\lambda]$ contains elements
of degree two (the two-cycles in $\tilde{Q}$).

And since $\tilde{Q}$ has two-cycles, the conclusion of Theorem
\ref{ref-3.1-3} does indeed not hold.  So hypothesis (5) cannot be deleted
from the statement of Theorem \ref{ref-3.1-3}.
\appendix
\section{Potentials associated to Del Pezzo surfaces}
\label{ref-A-15}
Let the notations be as in \S\ref{ref-4.2-10}. In this appendix we
will give an independent proof of Theorem \ref{ref-4.2.1-11}. We keep the setting of the
theorem although with trivial changes most of the proof applies to
arbitrary Fano varieties. Our arguments are related to those appearing
in recent papers by Minamoto and Mori \cite{minamoto,minamotomori}.
\begin{lemma}
\label{ref-A.1-16}
$B(\HH)$ is the derived tensor algebra of $B(\HH)_1$ over $B(\HH)_0=A(\EE)$.
\end{lemma}
\begin{proof}
For $k<0$ we have
\[
\Hom(\EE,\omega_Y^{-k}\otimes \EE)=\Ext^2(\EE,\omega_Y^{k+1}\otimes \EE)^\ast=0
\]
For $k\ge 0$ we have
\[
\Hom(\EE,\omega_Y^{-k}\otimes \EE)=R\Hom(\EE,\omega_Y^{-k}\otimes \EE)
\]
For $k,l\ge 0$ we have to prove that the canonical map
\[
\RHom(\EE,\omega_Y^{-k}\Lotimes \EE)\Lotimes_{A(\EE)}\RHom(\EE,\omega_Y^{-l}\otimes \EE)\r
\RHom(\EE,\omega_Y^{-k-l}\otimes \EE)
\]
is an isomorphism

One verifies that this is equivalent to showing that the composition map
\[
\RHom(\EE,\omega_Y^{-k}\Lotimes \EE)\Lotimes_{A(\EE)}\RHom(\omega_Y^{l}\otimes\EE,\EE)\r
\RHom(\omega_Y^{l}\otimes \EE,\omega_Y^{-k}\otimes \EE)
\]
is an isomorphism.

Since $\EE$ is a classical generator for $D^b(\coh(\Oscr_Y))$ it suffices
to prove the result
with $\omega_Y^{-k}\otimes \EE$ and $\omega_Y^{l}\otimes\EE$ replaced by
$\EE$. Then the statement is obvious.
\end{proof}
\begin{lemma}
\label{ref-A.2-17} Consider $A(\EE)\otimes A(\EE)$ as $A(\EE)$-bimodule with the outer
bimodule structure.
We have an equality of $A(\EE)$-bimodules 
\[
B(\HH)_1=\RHom_{A(\EE)^e}(A(\EE),A(\EE)\otimes A(\EE))[2]
\]
where on the righthandside the $A(\EE)$-bimodule structure comes from the surviving
inner bimodule structure on $A(\EE)\otimes A(\EE)$. 
In other words $B(\HH)_1$ is a shifted inverse dualizing complex of $A(\EE)$ (\cite[\S3.3]{Keller11}).
\end{lemma}
\begin{proof}
We have to compute
\[
\RHom_{\End(\EE)^e}(\End(\EE),\End(\EE)\otimes \End(\EE))
\]
Working on $Y\times Y$ we have
\[
\End(\EE)^e\cong \End(\EE\boxtimes \EE^\ast)
\]
It is easy to see that $\EE\boxtimes \EE^\ast$ is a classical generator for
$D^b(\coh(Y\times Y))$ (e.g.\ because it is derived from an exceptional collection, or
else invoke \cite{BondalVdB}). Under the ensuing derived equivalence
$D^b(\End(\EE\boxtimes \EE^\ast))\cong D^b(\coh(Y\times Y))$ it is easy
to see that $\End(\EE)$ corresponds to $\Oscr_Y$ and $\End(\EE)\otimes \End(\EE)$
to $\EE\boxtimes \EE^\ast$. 

We now have to compute
\[
\RHom_{Y\times Y}(\Oscr_\Delta,\EE\boxtimes \EE^\ast)
\]
where $\Delta$ is the diagonal. We use the well-known formula
\[
\uRHom_{\Oscr_{Y\times Y}}(\Oscr_\Delta,\Oscr_{Y\times Y})=\omega^{-1}_\Delta[-2]
\]
Thus
\begin{align*}
  \RHom_{Y\times Y}(\Oscr_\Delta,\EE\boxtimes \EE^\ast)&=
  R\Gamma(Y\times Y,\uRHom_{Y\times Y}(\Oscr_\Delta,\Oscr_{Y\times Y})
\Lotimes_{Y\times Y}\EE\boxtimes \EE^\ast)\\
  &=R\Gamma(Y\times Y,\omega^{-1}_\Delta\Lotimes_{Y\times Y}\EE\boxtimes \EE^\ast)[-2]\\
  &=\RHom_Y(\EE,\omega^{-1}_Y\otimes \EE)[-2]\qed
\end{align*}
\def\qed{}\end{proof}
\begin{lemma}
\label{ref-A.3-18}
Assume that $B$ is an $\NN$-graded ring such that $\dim B_i<\infty$ for all $i$.  Then
\[
\gldim B_0\le \gldim B
\]
If the category of graded finite length $B$-modules is $\Ext$-finite
and has a Serre functor given by $?(d)[n]$ for $d\neq 0$ then
\[
\gldim B_0< \gldim B=n
\]
\end{lemma}
\begin{proof}
Let $S,T$ be simple $B_0$-modules which we view as $B$-modules concentrated
in degree zero. The part of degree zero of a graded projective $B$-resolution
of $S$ is a projective $B_0$-resolution of $S$. This implies the first inequality. We
also get
\[
\Ext^i_{B_0}(S,T)=\Ext^i_B(S,T)_0
\]
Assume now that the category of graded finite length $B$-modules is $\Ext$-finite
and has a Serre functor given by $?(d)[n]$ for $d\neq 0$. Then it is standard
that $\gldim B=n$. For the second inequality we note that
$\Ext^n_B(S,T)_0=\Hom_B(T,S(d))^\ast_0=0$.
\end{proof}
\begin{proof}[Proof of Theorem \ref{ref-4.2.1-11}]
Since $A(\EE)$ is finite dimensional and has finite global dimension
the same is true for $A(\EE)^e$. Hence $A(\EE)$ is perfect as a bimodule 
and as such homologically smooth.

By Lemmas \ref{ref-A.1-16},\ref{ref-A.2-17} we obtain that $B(\HH)$ is a (graded)
$3$-Calabi-Yau completion of $A(\EE)$ (see \cite[\S4]{Keller11}). By
\cite[Thm 4.8]{Keller11} we obtain that $B(\HH)$ is 3-Calabi-Yau and
the proof shows that the Serre functor is given by $?(-1)$. By Lemma
\ref{ref-A.3-18} we then find that $\gldim A(\EE)\le 2$.  We may then use
\cite[Thm 6.10]{Keller11} to obtain that $B(\HH)$ is quasi-isomorphic to
a Ginzburg algebra derived from a super potential of degree one.
Since $B(\HH)$ is concentrated in degree zero it is also a Jacobian algebra. 

To prove that $B(\HH)$ is almost Azumaya we note that $B(\HH)$ is
the pushforward of an Azumaya algebra on the canonical bundle of $Y$. This
easily yields the desired result.
\end{proof}


\def\cprime{$'$} \def\cprime{$'$} \def\cprime{$'$} \def\cprime{$'$}
\ifx\undefined\bysame
\newcommand{\bysame}{\leavevmode\hbox to3em{\hrulefill}\,}
\fi

\end{document}